\documentclass[11pt, oneside]{article} 

\usepackage[utf8]{inputenc}
\usepackage[T1]{fontenc}
\usepackage{xcolor} 
\usepackage{color} 

\usepackage{amsmath}
\usepackage{amssymb} 
\usepackage{mathrsfs} 
\usepackage{bm} 
\usepackage{bbm} 

\usepackage{amsthm} 
\theoremstyle{plain}
\newtheorem{theorem}{Theorem} 

\newtheorem{condition}[theorem]{Condition}

\usepackage{algorithm2e} 
\SetKwInput{KwInputs}{Inputs}
\usepackage[noend]{algpseudocode}
\makeatletter
\def\BState{\State\hskip-\ALG@thistlm}
\makeatother

\usepackage{graphicx} 
\usepackage{subcaption} 

\usepackage{tikz} 
\usetikzlibrary{positioning ,shadows,matrix}

\usepackage[english]{babel} 

\usepackage{hyperref} 
\hypersetup{
    pdfpagemode={UseOutlines},
    bookmarksopen,
    pdfstartview={FitH},
    colorlinks,
    linkcolor={blue}, 
    citecolor={blue}, 
    urlcolor={blue} 
}

\usepackage[letterpaper,margin=1.25in]{geometry} 

\usepackage{enumitem} 

\def\iid{\overset{\textnormal{iid}}{\sim}} 

\makeatletter  
\@ifundefined{@currsizeindex} 
  {\RequirePackage{relsize}\let\dolarger\relsize} 
  {\def\dolarger#1{\larger[#1]}} 
\newcommand*\@@bigtimes[2]{\vphantom{\prod} 
  \vcenter{\hbox{\dolarger{4}$\m@th#1\mkern-2mu\times\mkern-2mu$}}} 
\newcommand*\bigtimes{\mathop{\mathpalette\@@bigtimes\relax}\displaylimits} 
\makeatother

\def\N{\mathbb{N}}\def\R{\mathbb{R}}\def\1{\mathbbm{1}}

\def\Fcal{\mathcal{F}}\def\Hcal{\mathcal{H}}\def\Rcal{\mathcal{R}}\def\Scal{\mathcal{S}}\def\Tcal{\mathcal{T}}\def\Wcal{\mathcal{W}}

\title{\bf Gaussian Process Methods for Covariate-Based Intensity Estimation}
\author{Patric Dolmeta and Matteo Giordano \\ \\ ESOMAS Department, University of Turin}
\date{} 

\begin{document}

\maketitle

\abstract{
We study nonparametric Bayesian inference for the intensity function of a covariate-driven point process.
We extend recent results from the literature, showing that a wide class of Gaussian priors, combined with flexible link functions, achieve minimax optimal posterior contraction rates. Our result includes widespread prior choices such as the popular Matérn processes, with the standard exponential (and sigmoid) link, and implies that the resulting methodologically attractive procedures optimally solve the statistical problem at hand, in the increasing domain asymptotics and under the common assumption in spatial statistics that the covariates are stationary and ergodic.
}

\bigskip

\noindent\textbf{Keywords.} Cox process; frequentist analysis of Bayesian procedures; Gaussian prior; minimax optimality

\tableofcontents

%

\section{Introduction} \label{sec:Intro}

In the statistical analysis of many point patterns, a central issue is to determine the influence of covariates on the point distribution. We refer to \cite{BCST12}, and references therein, for an overview of the problem and applications. The established framework to tackle this scenario is to model the point pattern via a Cox process  \cite{C55}, namely a `doubly-stochastic' point process with covariate-dependent intensity function, cf.~Section \ref{Sec:ObsModel} below. The resulting statistical problem is then to estimate the intensity from observations of the points and covariates.

	There is a vast literature on the parametric approach to this issue. For example the log-Gaussian Cox model \cite{MSW98} is widely used; further see the monograph \cite{D14}. Existing nonparametric frequentist  approaches largely rely on kernel-type estimators, e.g.~\cite{BCST12}. These methods were first shown to be asymptotically consistent in a seminal paper by Guan \cite{G08}, under the typical `large domain' framework of spatial statistics and the assumption that the covariates are stationary and ergodic.

	On the other hand, nonparametric Bayesian intensity estimation has so far been considered almost exclusively in models without covariates; see \cite{MSW98,AMM09,BSvZ15,KvZ15,DRRS17} for methodology and theoretical results, as well as for further references. The first frequentist asymptotic analysis of posterior distributions for covariate-driven point processes was developed in the recent paper by Giordano et al.~\cite{GKR23}, in an analogous increasing domain framework as the one considered in \cite{G08}. Among their results, they showed that truncated Gaussian wavelet priors, combined with suitable positive link functions, achieve minimax-optimal rates of posterior contraction towards the ground truth, in $L^1$-distance.

	In this paper, we build on the investigation of \cite{GKR23}, obtaining in the same setting optimal posterior contraction rates for a much wider class of (rescaled) Gaussian priors, including widespread choices such as the popular Matérn processes (e.g.~\cite[Section 11.4.4]{GvdV17}). See Theorem \ref{Theo:Main}. The proof of our result follows the general program set forth in \cite{GKR23}, based on the derivation of preliminary convergence properties in a covariate-dependent metric, cf.~\eqref{Eq:EmpDist}, to be combined with concentration inequalities for integral functional of stationary processes. To obtain the required uniformity over the prior support beyond the specific wavelet structure considered in \cite{GKR23}, we then employ techniques from the statistical theory of inverse problems, e.g.~\cite{GN20,GR22,N23}.

	Another novelty in our result is that, compared to \cite{GKR23}, we also significantly weaken the assumptions on the link function, which we only require to be locally Lipschitz (and bijective). Among the others, this allows to incorporate into the procedure the standard exponential link, as well as the sigmoid one, which had previously been employed for efficient posterior sampling in non-covariate-based intensity estimation \cite{AMM09}; see also \cite{KvZ15}. Thus, while theoretical, our result has also potential interesting methodological implications. We discuss these and other related research question in the concluding discussion in Section \ref{Sec:Discussion}.
	
%
%
%
%
%

\section{Posterior contraction rates for covariate-based intensity estimation}

\subsection{Observation model}
\label{Sec:ObsModel}

Throughout, we take $\R^D$ ($D\in\N$) to be the ambient space. Let the covariates be given by a 
$d$-dimensional ($d\in\N$) random field $Z := (Z(x), \ x\in\R^D)$, and let $(\Wcal_n)_{n\in\N}$ be a sequence of compact
sets in $\R^D$ satisfying $\Wcal_n\subseteq \Wcal_{n+1}$ for all $n$.

	For some $n\in\N$, the data consist in an observed realisation of the covariates over the window $\Wcal_n$, denoted by $Z^{(n)}:=(Z(x), \ x\in \Wcal_n)$, and of a random point process $N^{(n)}$ on $\Wcal_n$ which, conditionally given $Z^{(n)}$, is of inhomogeneous Poisson type with intensity
\begin{equation}
\label{Eq:LambdaRho}
	\lambda_\rho^{(n)}(x) := \rho(Z(x)), \qquad x\in\Wcal_n,
\end{equation}
 for some unknown 
 bounded  function $\rho : \R^d \to [0,\infty)$. Formally, we may write
\begin{equation}
\label{Eq:PointProc}
	N^{(n)} \overset{d}{=} \{X_1,\dots, X_{N_n}\}, 
	\quad N_n|Z^{(n)}\sim\textnormal{Po}(\Lambda^{(n)}_\rho),
	\quad
	X_i|Z^{(n)} \iid \frac{\lambda^{(n)}_\rho(x)dx}{\Lambda^{(n)}_\rho},
\end{equation}
with $\Lambda^{(n)}_\rho := \int_{\Wcal_n}\lambda^{(n)}_\rho(x) dx$. In other words, $N^{(n)}$ is a Cox process on $\Wcal_n$ directed by the random measure $\lambda^{(n)}_\rho(x)dx$ \cite{C55}. The statistical problem is to estimate non-parametrically $\rho$ from data $(Z^{(n)},N^{(n)})$.

	We will denote by $P^{(n)}_{\rho}$ the joint law of $(N^{(n)},Z^{(n)})$. By standard theory for Poisson processes, e.g.~\cite{K98}, $P^{(n)}_{\rho}$ is absolutely continuous with respect to the law $P^{(n)}_1$ corresponding to the homogeneous case, with log-likelihood equal to
\begin{equation}
\label{Eq:Likelihood}
	l_n(\rho):=\log\frac{dP^{(n)}_\rho}{dP^{(n)}_1}(D^{(n)})=
	\int_{\Wcal_n} \log(\rho(Z(x))) d N^{(n)}(x) 
	-\int_{\Wcal_n}\rho(Z(x))dx,
\end{equation}
cf.~\cite[Theorem 1.3]{K98}. We will write $E_\rho^{(n)}$ for the expectation with respect to $P_\rho^{(n)}$.

%
%
%

\subsection{Assumptions on the covariates}

We consider the important case where the covariates are given by a Gaussian process, according to the following assumption. Under the latter, the point process \eqref{Eq:PointProc} is then a non-parametric generalisation of the celebrated log-Gaussian Cox process \cite{MSW98}, wherein $\rho$ in \eqref{Eq:LambdaRho} takes the parametric form $\rho(z) = \exp(\beta^Tz)$ for some $\beta\in\R^d$.  Gaussian processes are among the most ubiquitously employed models for spatial random fields in statistical applications, e.g.~\cite[Section 2.3]{C15}.

\begin{condition}\label{Cond:GaussCov}
Let $\tilde Z^{(h)}:=(\tilde Z^{(h)}(x), \ x\in\R^D)$, for $h=1,\dots,d,$ be independent, almost surely locally bounded, centred and stationary Gaussian processes with integrable covariance functions. Further assume without loss of generality that $\textnormal{Var}[Z^{(h)}(x)]=1$ for all $h$ and $x$. Let $Z$ be given by 
\begin{equation}
\label{Eq:PhiZ}
	Z(x) := [\phi(\tilde Z^{(1)}(x)), \dots, \phi(\tilde Z^{(d)}(x))], \qquad x\in\R^D,
\end{equation}
where $\phi$ is the standard normal cumulative distribution function.
\end{condition}

Under Condition \ref{Cond:GaussCov}, $\tilde Z(x):=(\tilde Z^{(1)}(x),\dots,\tilde Z^{(d)}(x))\sim N_d(0,1)$ for all $x\in\R^D$, which amounts to the common practice of standardising the covariates, cf.~\cite[Section 3.2.1]{G08}. The component-wise application of $\phi$ in \eqref{Eq:PhiZ} is without loss of generality in view of its invertibility. It can be thought of as a convenient `pre-processing' step that implies, in particular, that $Z$ takes values in the compact set $[0,1]^d$. The integrability of the covariances is a mild requirement that is satisfied under the (often realistic, e.g.~\cite[Sec.~2.3]{C15}) assumption that the correlation between covariates at distant locations decays sufficiently fast. This condition implies that $Z$ is ergodic, which, similarly to \cite{G08}, crucially enables in our analysis a sufficient `accumulation of information' from distinct points in the observation window with similar covariate values, allowing for consistent inference in the increasing domain asymptotics.

%
%
%

\subsection{Prior specification}

Gaussian processes are among the most universal models for prior distributions in function spaces, and have been successfully employed in related nonparametric Bayesian intensity estimation problems; see \cite{MSW98,AMM09,KvZ15} and references therein.

	In the present setting, under Condition \ref{Cond:GaussCov}, we place on the unknown function $\rho :[0,1]^d \to[0,\infty)$ in \eqref{Eq:LambdaRho} an ($n$-dependent) log-Gaussian prior $\Pi_n$, constructed starting from a base distribution that we require to satisfy the following mild regularity condition.

\begin{condition}\label{Cond:GPCondition}
	Let $\tilde \Pi$ be a centred Gaussian probability
measure on the Banach space $C([0,1]^d)$ that is supported on a linear subspace of $C^1([0,1]^d)$. Further assume that, for some $\alpha>0$, the reproducing kernel Hilbert space (RKHS) $\tilde \Hcal$ of $\tilde \Pi$ is equal to the Sobolev space $H^{\alpha}([0,1]^d)$, with RKHS norm satisfying $\|\cdot\|_\Hcal\simeq \|\cdot\|_{H^\alpha}$. 
\end{condition}

	Condition \ref{Cond:GPCondition} is satisfied by a wide variety of Gaussian processes, including by stationary ones with polynomially-tailed spectral measures (e.g.~the popular Matérn processes, cf.~\cite[Section 11.4.4]{GvdV17}), as well as by series priors defined on bases spanning the scale of Sobolev spaces (e.g.~\cite[Section 11.4.3]{GvdV17}).

	Given $\tilde \Pi$ satisfying Condition \ref{Cond:GPCondition}, we then construct the prior $\Pi_n$ as the law of the random function
\begin{equation}
\label{Eq:LogGP}
	R_n(z) := \exp\left( n^{-d/(4\alpha+2d)} \tilde W(z)\right), \qquad z\in[0,1]^d,
	\qquad \tilde W\sim\tilde\Pi.
\end{equation}
The introduction of the exponential link in \eqref{Eq:LogGP} is motivated by the positivity constraint on $\rho$. This particular choice was made for concreteness, and we note that any other locally Lipschitz and bijective link function may be used in the result to follow, including for example the sigmoid one used in \cite{AMM09,KvZ15}. The $n$-dependent rescaling of the base Gaussian prior is a common technique in the statistical theory for inverse problems, e.g.~\cite{GN20,GR22,N23}. In the proof, it implies a bound for the norm of the draws from $\Pi$, which is key to obtaining certain uniform concentration inequalities, going beyond the specific wavelet structure considered in \cite{GKR23}. Non-asymptotically, the rescaling amounts to a simple adjustment of the covariance function of $\tilde\Pi$, and does not require further tuning in practice.

%
%
%

\subsection{Main result}

In our main result, we characterise the speed of asymptotic concentration of the posterior distributions $\Pi_n(\cdot|N^{(n)},Z^{(n)})$ arising from the log-Gaussian priors \eqref{Eq:LogGP}, under the frequentist assumption that the data $(N^{(n)},Z^{(n)})$ have been generated by some fixed ground truth $\rho_0$. We work in the (widely adopted in spatial statistics) increasing domain asymptotics, namely $\textnormal{vol}(\Wcal_n)\to\infty$ as $n\to\infty$, which reflects the common situation where a single realisation of the points and covariates is observed over a large window. Similarly to \cite{G08}, our proof requires $\Wcal_n$ to grow uniformly in all directions, and for concreteness we restrict to square domains
\begin{align}
\label{Eq:SpatialWn}
	\Wcal_n=\left[-\frac{1}{2}n^{1/D},\frac{1}{2}n^{1/D}\right]^D,
\end{align}
whence $\textnormal{vol}=n$. Extensions to other domains containing the above sets is straightforward.

\begin{theorem}\label{Theo:Main}
Let $\rho_0 = \exp\circ w_0$, for some $w_0\in C^\beta([0,1]^d)\cap H^\beta([0,1]^d)$ and some $\beta>\min(1,d/2)$. Consider data $(N^{(n)},Z^{(n)})\sim P^{(n)}_{\rho_0}$ arising from the observation model \eqref{Eq:PointProc} with $\rho=\rho_0$, $\Wcal_n$ as in \eqref{Eq:SpatialWn} and with $Z$ a random field satisfying Condtion \ref{Cond:GaussCov}. Let the prior $\Pi_n$ be as in \eqref{Eq:LogGP}, with $\tilde\Pi$ satisfying Condition \ref{Cond:GPCondition} with $\alpha=\beta$. Then, for $M>0$ large enough, as $n\to\infty$,
$$
	E_{\rho_0}^{(n)}
	\Bigg[\Pi_n\Big(\rho : \|\rho - \rho_0\|_{L^1} > M n^{-\beta/(2\beta+d)}
	\Big| N^{(n)},Z^{(n)}\Big)\Bigg]
	\to 0.
$$
\end{theorem}

	Theorem \ref{Theo:Main} asserts that the posterior distribution $\Pi_n(\cdot|N^{(n)},Z^{(n)})$ asymptotically concentrates its mass over small neighbourhoods of the ground truth, whose $L^1$-radius shrinks at rate $n^{-\beta/(2\beta+d)}$. The latter is known to be minimax for $\beta$-smooth intensities, e.g.~\cite{K98}, implying the optimality of Theorem \ref{Theo:Main}.

	The proof of the result, presented in Section \ref{Sec:Proof}, is based on recent work by Giordano et al.~\cite{GKR23}, who developed a general program to derive posterior contraction rates in covariate-based intensity estimation. In particular, the argument also implies that $\Pi(\cdot|N^{(n)},Z^{(n)})$ concentrates at the same speed also in the `empirical' (i.e.~covariate-dependent) $L^1$-metric
\begin{equation}
\label{Eq:EmpDist}
	\frac{1}{n}\|\lambda_\rho^{(n)} - \lambda_{\rho_0}^{(n)}\|_{L^1(\Wcal_n)}
	=\frac{1}{|\Wcal_n|}\int_{\Wcal_n}|\rho(Z(x))-\rho_0(Z(x))|dx,
\end{equation}
cf.~\eqref{Eq:Intermediate} below. This implies that the proposed nonparametric Bayesian methodology also optimally solves the `prediction' problem of inferring the overall spatial intensity function $\lambda_\rho^{(n)}$ in \eqref{Eq:LambdaRho}.
	
%
%
%
%
%
%

\section{Concluding discussion}
\label{Sec:Discussion}

In this article, we have considered Gaussian process methods for estimating the intensity function of a covariate-driven point process. In a realistic increasing domain framework with ergodic covariates, our main result, Theorem \ref{Theo:Main}, shows that the posterior distributions arising from a large class of (rescaled) Gaussian priors concentrates at optimal rate around the ground truth. This is based on recent work by \cite{GKR23}, and extends a result obtained in the latter reference for the specific case of truncated Gaussian wavelet priors.

	We conclude mentioning some related interesting research questions. Firstly, Theorem \ref{Theo:Main} is non-adaptive in that the prior specification under which optimal rates are obtained requires knowledge of the regularity $\beta$ of the ground truth. Extensions to adaptive version of the rescaled priors \eqref{Eq:LogGP} in the present setting is of primary interest, but is known to be an involved issue in the broader literature, cf.~\cite{GN20,GR22}.

	Implementation of posterior inference with the considered Gaussian priors (as well as the other priors studied in \cite{GKR23}) is also of great importance. In ongoing work, we are exploring extensions to the problem at hand of the existing efficient posterior sampling methodologies for non-covariate-based Bayesian intensity estimation \cite{AMM09}.

	Lastly, beyond the large domain framework considered in the present article, `infill asymptotics' are also of interest in spatial statistics. For covariate-based intensity estimation, this setting can be formulated as having access to repeated (and possibly independent) realisations of the covariates and the points over a fixed observation window. We refer to \cite{GKR23} for possible extensions of the proof techniques employed here to such situations. This is beyond the scope of the present article, but represents an interesting avenue for future research.

\section{Proof}
\label{Sec:Proof}

We employ techniques from \cite{GKR23}, to which we refer for the necessary background. We begin proving posterior contraction in the empirical distance \eqref{Eq:EmpDist}. This is then combined with an exponential concentration inequality for 
$n^{-1}\|\lambda_\rho^{(n)} - \lambda_{\rho_0}^{(n)}\|_{L^1(\Wcal_n)}$ around its `ergodic average' $\|\rho - \rho_0\|_{L^1}$, which allows to extend the obtained rates to the latter metric, concluding the proof.

	Set $\varepsilon_n:= n^{-\beta/(2\beta+d)}$. For $\tilde\Pi$ satisfying Condition \ref{Cond:GPCondition}, let $\tilde W\sim \tilde\Pi$ and denote by $\tilde\Pi_n$ the law of the random function $(\sqrt n \varepsilon_n)^{-1}\tilde W = n^{-d/(4\alpha+2d)}\tilde W$.
Using standard techniques from the posterior contraction rate theory for (rescaled) Gaussian priors, e.g.~\cite{GN20,GR22,N23}, for all $K_1>0$ there exist a sufficiently large constant $M_1>0$ and a set 
$$
	\Fcal_n\subset \left\{ f\in C^1([0,1]^d) : \|f\|_{C^1}\le M_1\right\}
$$
such that, for all $n$ large enough,
$$
	\tilde \Pi_n (\Fcal_n^c)\le e^{-K_1 n\varepsilon_n^2}; \qquad
	\log N(\varepsilon_n;\Fcal_n,\|\cdot\|_{L^\infty})\lesssim n\varepsilon_n^2.
$$
See e.g.~Lemma 5 in \cite{GR22} and its proof. Set $\Rcal_n:=\left\{\exp\circ f,\ f\in \Fcal_n\right\}$. Then, by construction of $\Pi_n$, cf.~\eqref{Eq:LogGP}, the first inequality in the last display implies that $\Pi_n (\Rcal_n^c)\le e^{-K_1 n\varepsilon_n^2}$. Further, since $\Fcal_n$ is contained in a sup-norm ball, $\| \exp\circ f_1 - \exp\circ f_2\|_{L^\infty}\lesssim \|f_1 -  f_2\|_{L^\infty}$ for all $f_1,f_2\in\Fcal_n$. It follows that
$
	\log N(\varepsilon_n;\Rcal_n,\|\cdot\|_{L^\infty})
	 \lesssim \log N(\varepsilon_n;\Fcal_n,\|\cdot\|_{L^\infty})\lesssim n\varepsilon_n^2.
$
Lastly, using Corollary 2.6.18 of \cite{GN16}, since $w_0\in H^\beta([0,1])^d=\Hcal$ by assumption, we have
\begin{align*}
	\tilde\Pi_n(w:\|w - w_0\|_{L^\infty} \le \varepsilon_n)
	\ge e^{-\frac{1}{2}\|w_0\|_{\Hcal}^2n\varepsilon_n^2}\tilde\Pi(w:\|w\|_\infty\le \sqrt n\varepsilon_n^2),
\end{align*}
which, by the metric entropy estimate in Theorem 4.3.36 of \cite{GN16} and the centred small ball inequality in Theorem 1.2 of \cite{LL99}, is lower bounded by
$
	e^{-\frac{1}{2}\|w_0\|_{\Hcal}^2n\varepsilon_n^2}
	e^{-c_1n\varepsilon_n^2}
	=e^{-c_2n\varepsilon_n^2}
$
for some $c_1,c_2>0$. Arguing as in the proof of Theorem 3.2 in \cite{GKR23}, this implies that for all sufficiently large $M_2>0$, as $n\to\infty$,
\begin{equation}
\label{Eq:Intermediate}
	E_{\rho_0}^{(n)}
	\left[\Pi(\Scal_n| N^{(n)},Z^{(n)})\right]
	\to 1,
\end{equation}
where
$
	\Scal_n := \left\{\rho\in\Rcal_n 
	: n^{-1}\|\lambda_\rho^{(n)} - \lambda_{\rho_0}^{(n)}\|_{L^1(\Wcal_n)} \le  M_2 
	\varepsilon_n\right\}.
$
Now for $M>0$ to be chosen below, set 
$
	\Tcal_n :=\left\{\rho\in C^1([0,1]^d):\|\rho - \rho_0\|_{L^1} \le M\varepsilon_n\right\}.
$
The proof is then concluded by showing that $E_{\rho_0}^{(n)} [\Pi( \Tcal_n^c | N^{(n)}, Z^{(n)}) ] \to0$ as $n\to\infty$. To do so, note that by \eqref{Eq:Intermediate},
\begin{align*}
 	\Pi( \Tcal_n^c | N^{(n)}, Z^{(n)}) 
	&= \Pi( \Tcal^c_n \cap \Scal_n | N^{(n)}, Z^{(n)}) + o_{P^{(n)}_{\rho_0}}(1)\\ 
	&= \frac{ \int_{\Tcal^c_n \cap \Scal_n} e^{l_n(\rho) - l_n(\rho_0)}d\Pi(\rho) }
	{ \int_{\mathcal C([0,1]^d)} e^{l_n(\rho) - l_n(\rho_0)  } d\Pi(\rho)} 
	+ o_{P^{(n)}_{\rho_0}}(1).
\end{align*}
Denote by $D_n$ the denominator in the previous display. The proof of Theorem 3.2 in \cite{GKR23} shows that $P_{\rho_0}^{(n)}( D_n \leq e^{ - K_1 n \varepsilon_n^2 } ) = o(1)$ as $n\to\infty$ for some constant $K_1>0$, so that by Fubini's theorem,
\begin{align*}
 	&E_{\rho_0}^{(n)} \left[\Pi( \Tcal_n^c | N^{(n)}, Z^{(n)})\right] \\
	&\ \leq e^{  K_1 n \varepsilon_n^2 }  
 	\int_{\Tcal_n^c \cap \{\rho\in\mathcal R_n:\|\rho\|_{C^1}\le M_1\}} 
	\Pr (\|\lambda^{(n)}_\rho - \lambda^{(n)}_{\rho_0}
	\|_{L^1(\Wcal_n)} \le M_2n\varepsilon_n ) d\Pi(\rho)+ o(1).
\end{align*}
Fix any $\rho\in \Tcal_n^c \cap \{\rho\in\mathcal R_n:\|\rho\|_{C^1}\le M_1\}$. Then, if $\|\lambda^{(n)}_\rho - \lambda^{(n)}_{\rho_0}\|_{L^1(\Wcal_n)} \le M_2n\varepsilon_n $, necessarily 
$$
	\Delta^{(n)}(\rho) := \| \rho - \rho_0\|_{L^1} 
	- \frac{1}{n}\|\lambda^{(n)}_\rho - \lambda^{(n)}_{\rho_0}
	\|_{L^1(\Wcal_n)} > (M-M_2)\varepsilon_n
	\geq \frac{M_2}{2}\varepsilon_n  
$$ 
upon taking $M > 2M_2$. For all such $M$, it follows that the expectation of interest is upper bounded by
\begin{align*}
	e^{  K_1 n \varepsilon_n^2 }  \int_{\Tcal_n^c \cap \{\rho\in\mathcal R_n:\|\rho\|_{C^1}\le M_1\}}
	 \Pr (\Delta^{(n)}(\rho) > M_2\varepsilon_n/2 )d\Pi(\rho)+ o(1)  .  
\end{align*}
The concentration inequality in Proposition D.1 of \cite{GKR23}, applied with $f := |\rho - \rho_0|- n^{-1}\|\lambda^{(n)}_\rho - \lambda^{(n)}_{\rho_0}\|_{L^1(\Wcal_n)}$, for $\rho\in C^1([0,1]^d)$, whose (weak) gradient satisfies $ \| \nabla f  \|_{L^\infty}\le \| \rho  \|_{C^1}+\|  \rho_0  \|_{C^1}\le 2M_1$ now gives that
$$
 	\sup_{\rho\in \Tcal_n^c \cap \{\rho\in\mathcal R_n:\|\rho\|_{C^1}\le M_1\}}
	\Pr (\Delta^{(n)}(\rho) > M_2\varepsilon_n/2 )\leq e^{- K_2 (M_2)^2 n\varepsilon_n^2 }
$$
for some $K_2>0$. The claim follows taking $M_2>0$ large enough and combining the last two displays.\qed

\paragraph{acknowledgement}
M.G.~has been partially supported by MUR, PRIN project 2022CLTYP4.

\bibliographystyle{acm}

\bibliography{PointProcRefIWFOS25}

\end{document}